\documentclass[12pt,reqno]{amsart}
\usepackage{amsmath,amsfonts,amsthm,amssymb,amsxtra}
\usepackage{graphicx}
\newcommand\version{November 8, 2007} 

\usepackage{amsxtra, amssymb, mathrsfs}

\newtheorem{theorem}{Theorem}[section]
\newtheorem{proposition}[theorem]{Proposition}
\newtheorem{lemma}[theorem]{Lemma}

\theoremstyle{definition}

\newtheorem{example}[theorem]{Example}

\theoremstyle{remark}

\newtheorem{remark}[theorem]{Remark}

\numberwithin{equation}{section}

\newcommand{\A}{\mathcal{A}}

\newcommand{\C}{\mathbb{C}}

\newcommand{\E}{\mathcal{E}}

\renewcommand{\H}{H}

\newcommand{\loc}{{\rm loc}}
\newcommand{\N}{\mathbb{N}}
\newcommand{\Neumann}{\mathcal{N}}
\newcommand{\qf}{\mathfrak{a}}
\newcommand{\R}{\mathbb{R}}

\newcommand{\Z}{\mathbb{Z}}

\DeclareMathOperator{\im}{Im}
\DeclareMathOperator{\re}{Re}
\DeclareMathOperator{\spec}{spec}

\title[Hardy inequalities on metric trees --- \version]{Remarks about Hardy inequalities on metric trees}

\author{Tomas Ekholm}
\address{Tomas Ekholm, Centre for Mathematical Sciences, Lund
  University, Box 118, 22100 Lund, Sweden} 
\email{tomase@maths.lth.se}

\author{Rupert L. Frank}
\address{Rupert L. Frank, Department of Mathematics, Fine Hall,
  Princeton University, Princeton, NJ 08544, USA} 
\email{rlfrank@math.princeton.edu}

\author{Hynek Kova\v r\'{\i}k}
\address{Hynek Kova\v r\'{\i}k, Department of Mathematics, Stuttgart
  University, Pfaffen\-wald\-ring 57, 70569 Stuttgart, Germany} 
\email{hynek.kovarik@mathematik.uni-stuttgart.de} 

\begin{document}

\begin{abstract}
We find sharp conditions on the growth of a rooted regular metric tree such that the Neumann Laplacian on the tree satisfies a Hardy inequality. In particular, we consider homogeneous metric trees. Moreover, we show that a non-trivial Aharonov-Bohm magnetic field leads to a Hardy inequality on a loop graph.
\end{abstract}

\keywords{Laplace operator, metric tree, Hardy inequality}

\thanks{\copyright\, 2007 by the authors. This paper may be  
reproduced, in its entirety, for non-commercial purposes.}

\maketitle


\section{Introduction}

Let $\Gamma$ be a rooted metric tree of infinite height with root $o$.
We are interested in Hardy inequalities of the form
\begin{align} \label{model}
\int_{\Gamma} \psi(|x|) |u(x)|^2 \, dx
\leq C(\Gamma,\psi) \int_{\Gamma} |u'(x)|^2\, dx
\end{align}
where $\psi > 0$ is the so-called Hardy weight and $C(\Gamma,\psi)$ is a positive
constant which might depend on $\Gamma$ and $\psi$, but which is independent of
$u$. Evans, Harris and Pick \cite{EHP} found a necessary and sufficient condition such that \eqref{model} 
holds for all functions $u$ on $\Gamma$ such that $u(o)=0$ and such that the integral on the right hand side is finite. They consider even the case of non-symmetric Hardy weights. However, due to the existence of harmonic functions with finite Dirichlet integral, the Hardy weights have to decay rather fast. This led Naimark and Solomyak \cite{NS2} to the study of \eqref{model} for functions in $\{u \in C_0^\infty(\Gamma):\ u(o)= 0\}$ (and its closure with respect to the Dirichlet integral). For regular trees, see Subsection \ref{prelim} for the definition, they gave a complete characterization of the validity of \eqref{model} on that class of functions. Note that in both papers the condition $u(o)=0$ was imposed. It is clear that without this assumption inequality \eqref{model} cannot hold for all metric trees. To see this, it suffices to consider $\Gamma=\R_+$ as an example. 

Our first remark in this paper is that if the tree is regular and grows sufficiently fast, then there are weights $\psi$ such that \eqref{model} hold for all $u\in C_0^\infty(\Gamma\cup\{o\})$ (without the condition $u(o)=0$). Following the approach in \cite{NS2} we can give a complete characterization of admissible weights and obtain two-sided estimates on the sharp constant $C(\psi,\Gamma)$ in \eqref{model}. As in the Euclidean case, Hardy weights typically decay like $|x|^{-2}$ as $|x|:= \text{dist}(o,x)\to\infty$. The growth of a tree is reflected by its branching function 
\begin{equation}\label{eq:branching}
g_0(t) := \#\, \{x\in\Gamma: |x|=t\},\quad t\in\R_+\, ,
\end{equation}
and our condition for the validity of \eqref{model} is that the \emph{reduced height} of 
$\Gamma$ is finite, i.e.,
\begin{align}\label{reduced_height}
\int_0^{\infty}\, \frac{dt}{g_0(t)} < \infty\, .
\end{align}
These trees are usually called \emph{transient}. If $\Gamma$ has \emph{global dimension} $d$, that is,
\begin{align}\label{eq:dimension}
0 <  \inf_{t\geq0} \frac{g_0(t)}{(1+t)^{d-1}}\, \leq \,
    \sup_{t\geq0} \frac{g_0(t)}{(1+t)^{d-1}} <\infty\, ,
\end{align}
then \eqref{reduced_height} holds iff $d>2$.

Our second remark concerns tree with branching function growing faster than any power. As an example of this situation we study \emph{homogeneous trees}, i.e., regular trees where all the edges have the same length and all the vertices have the same branching number. In this case, the Laplacian is positive definite and one may ask whether \eqref{model} is still true if we subtract $\|u\|^2$ times the bottom of the spectrum from the right hand side. We prove that this is indeed the case, and give again a complete characterization of all admissible weights. Again the condition $u(o)=0$ is \emph{not} needed.

Our third remark in this paper points out another mechanism that induces a Hardy inequality. This time we consider not a tree, but a halfline with a loop attached. Clearly, for the Laplacian on such a graph no Hardy inequality is valid. We prove that an Aharonov-Bohm magnetic field with non-integer flux through the loop does give rise to a Hardy inequality. This is reminiscent of the Hardy inequality in the two-dimensional punctured plane with an Aharonov-Bohm magnetic field \cite{LW}.

\subsection*{Acknowledgements}

The authors are grateful to Timo Weidl for several useful discussions, and to the organizers of the workshop `Analysis on Graphs' at the Isaac Newton Institute in Cambridge for their kind invitation. This
work has been supported by FCT grant SFRH/BPD/ 23820/2005 (T.E.) and
DAAD grant D/06/49117 (R.F.). Partial support by the ESF programme
SPECT (T.E. and H.K.) and the DAAD-STINT PPP programme (R.F.) is
gratefully acknowledged.
 

\section{Hardy inequalities on regular transient trees}

\subsection{Preliminaries}
\label{prelim}
In this subsection we would like to recall some basic definitions about trees and fix our notation. We refer, e.g., to \cite{NS1,NS2,S} for details. Let $\Gamma$ be a rooted metric tree with root $o$. We denote by $|x|$ 
the unique distance between a point $x \in\Gamma$ and the root $o$ and always assume that $\Gamma$ is of infinite \emph{height}, i.e., $\sup_{x\in\Gamma}|x|=\infty$. The branching number $b(x)$ of a vertex
$x$ is defined as the number of edges emanating from $x$. We assume
the natural conditions that $b(x) > 1$ for any vertex $x \neq o$ and
that $b(o)=1$.

We will assume that $\Gamma$ is \emph{regular}, i.e., all the vertices
of the same distance to the root have equal branching numbers and all
the edges emanating from these vertices have equal length.

Let $x$ be a vertex such that there are $k+1$ vertices on the unique path
between $o$ and $x$ including the endpoints. We denote by $t_k$ the
distance $|x|$ and by $b_k$ the branching number of $x$. Moreover, we
put $t_0:=0$ and $b_0:=1$. Note that $t_k$ and $b_k$ are only
well-defined for regular trees and that these numbers, in the regular
case, uniquely determine the tree. 

Let the (first) branching function $g_0$ be defined by \eqref{eq:branching}, or equivalently,
\begin{align*}
  g_0(t) = b_0\, b_1\cdots b_k, 
  \quad \text{if} \ t_k < t \leq t_{k+1},
  \quad k\in \N\, .
\end{align*}
Note that $g_0$ is a non-decreasing function and that $g_0(t)$
coincides with the number of points $x \in \Gamma$ such that
$|x|=t$. The rate of growth of $g_0$ reflects the rate of growth of
the tree $\Gamma$. More precisely, $g_0$ measures how the surface of
the `ball' $\{x\in\Gamma : |x| < t\}$ grows with $t$. We shall assume that this growth is fast enough in the sense that \eqref{reduced_height} holds. Such trees are called \emph{transient} or of \emph{finite reduced height}.

The Sobolev space $\H^1(\Gamma)$ consists of all continuous functions
$u$ such that $u \in \H^1(e)$ on each edge $e$ of $\Gamma$
and 
\begin{align*}
\int_{\Gamma} \left(|u'(x)|^2 + |u(x)|^2\right) \, dx <
\infty\, ,
\end{align*}
and the Neumann Laplacian $-\Delta_\Neumann$ in $L_2(\Gamma)$ is defined through the quadratic form
\begin{align}\label{eq:kinetic}
  \int_{\Gamma} |u'(x)|^2 \, dx, \quad u \in \H^1(\Gamma)\,.
\end{align}
Functions in its domain satisfy a Neumann boundary condition at the root and Kirchhoff matching conditions at the vertices $x \neq o$.

\subsection{Hardy's inequality on transient trees} 

Let $\Gamma_o:=\Gamma\cup\{o\}$ and denote by $C_0^\infty(\Gamma_o)$ the class of infinitely smooth functions on $\Gamma_o$ the support of which is a bounded subset of $\Gamma_o$. We emphasize that functions in this class do not necessarily vanish at the root. We shall prove

\begin{theorem}\label{hardy_Mazya}
Let $\Gamma$ be a regular metric tree of infinite height and finite reduced height, and let $\psi$ be a measurable, non-negative function on $\R_+$. Then the Hardy inequality
\begin{equation} \label{eq:hardy_Mazya}
\int_{\Gamma} \psi(|x|) |u(x)|^2 \, dx
\leq C(\Gamma,\psi) \int_{\Gamma} |u'(x)|^2\, dx,
\quad u \in C_0^\infty(\Gamma_o),
\end{equation}
is valid if and only if
\begin{equation}\label{eq:hardycond_Mazya}
M(\Gamma,\psi) := \sup_{t>0} \left( \int_0^t \psi(s) g_0(s)\,ds \right) \left(
  \int_t^\infty \frac{ds}{g_0(s)} \right) < \infty.
\end{equation}
Moreover, the sharp constant in \eqref{eq:hardy_Mazya} satisfies 
$$
M(\Gamma,\psi) \leq C(\Gamma, \psi) \leq 4 M(\Gamma,\psi).
$$
\end{theorem}

It follows from \eqref{eq:hardy_Mazya} that the closure of $C_0^\infty(\Gamma_o)$ with respect to \eqref{eq:kinetic} is a function space and that \eqref{eq:hardy_Mazya} holds for all functions from this closure. In particular, it holds for all $u\in H^1(\Gamma)$. Let us give a simple

\begin{example}
Assume that $\Gamma$ has global dimension $d>2$ (see \eqref{eq:dimension}). Then one has
\begin{equation}\label{eq:hardydim}
\int_{\Gamma} \frac{|u(x)|^2}{(1+|x|)^2} \, dx
\leq C(\Gamma) \int_{\Gamma} |u'(x)|^2\, dx,
\quad u \in C_0^\infty(\Gamma_o)\, ,
\end{equation}
and the decay exponent of the Hardy weight cannot be improved.
\end{example}

\begin{remark}
  It is easy to see (and follows from the proof of Theorem \ref{hardy_Mazya}) that \emph{no} non-trivial Hardy inequality can hold on recurrent trees, i.e., trees for which the integral in \eqref{reduced_height} is infinite. This dichotomy is reminiscent of the Euclidean situation, where a Hardy inequality holds only in three and higher dimensions. Another manifestation of this fact is the validity of a Cwikel-Lieb-Rozenblum inequality for the number of negative eigenvalues of a Schr\"odinger operator on a tree, see \cite{EFK}. On a recurrent tree a Schr\"odinger operator has weakly coupled bound states, see \cite{K}.
\end{remark}


\subsection{Hardy's inequality on homogeneous trees}

It is a consequence of \eqref{eq:hardy_Mazya} (with $\psi\equiv 1$) that the Neumann Laplacian $-\Delta_\Neumann$ is positive definite iff
\begin{equation}\label{eq:hardycond_Mazya1}
M(\Gamma,1) = \sup_{t>0} \left( \int_0^t g_0(s)\,ds \right) \left(
  \int_t^\infty \frac{ds}{g_0(s)} \right) < \infty\, ,
\end{equation}
and that one has
$$
(4 M(\Gamma,1))^{-1} \leq \inf\spec(-\Delta_\Neumann) \leq M(\Gamma,1)^{-1} \, .
$$
Note that condition \eqref{eq:hardycond_Mazya1} coincides with the corresponding condition for the Dirichlet Laplacian. This is a priori not clear since the spectrum need not be purely continuous.

In the positive definite case it is natural to ask whether the Hardy inequality \eqref{eq:hardy_Mazya} continues to hold if we normalize the right hand side such that its spectrum starts at zero. In other words, can one replace 
$ \int |u'|^2\, dx$ by $ \int |u'|^2\, dx - \inf\spec(-\Delta_\Neumann) \int |u|^2\,dx$? We answer this question affirmatively for a special class of trees.
Recall that a regular tree $\Gamma$ is called \emph{homogeneous} if all the edges have the same length and if all the vertices have the same branching number $b > 1$. By scaling we can assume without loss of generality that
the edge length equals $1$. The branching function has then the form
\begin{align*}
g_0(t) = b^j, \quad j < t \leq j + 1, \quad j\in \N_0\, .
\end{align*}
A homogeneous tree clearly satisfies \eqref{reduced_height} and \eqref{eq:hardycond_Mazya1}. It turns out that the bottom of the spectrum of $-\Delta_\Neumann $ can be calculated explicitly (see \cite{SS}), namely,
\begin{align*}
\inf\spec(-\Delta_\Neumann) = \lambda_b := \left( \arccos \frac 1{R_b}\right)^2, \quad R_b := \frac
       {b^{\frac 12} + b^{-\frac 12}}{2}\, .
\end{align*}

We shall prove

\begin{theorem}\label{Thm:hardy_homo}
Let $\Gamma$ be a homogeneous metric tree with a branching number
$b\geq 2$ and edge length $1$, and let $\psi$ be a measurable, non-negative function on $\R_+$.
Then the Hardy inequality
\begin{align} \label{eq:hardy_homo}
\int_{\Gamma} \psi(|x|) |u(x)|^2 \, dx
\leq C(b,\psi) \int_{\Gamma} \left( |u'(x)|^2 -
\lambda_b|u(x)|^2\right)\, dx, \quad u \in C_0^\infty(\Gamma_o),
\end{align}
holds with some constant $C(b,\psi)$ if and only if
\begin{equation} \label{condition-homo}
\sup_{r>0}\, (1 + r)^{-1} \int_0^r \psi(t) (1+t)^2 \, dt < \infty\,.
\end{equation}
\end{theorem}

{F}or two-sided estimates on the sharp constant $C(b,\psi)$ in terms of \eqref{condition-homo} we refer to the proof below. In analogy to \eqref{eq:hardydim} we record
\begin{align*}
\int_{\Gamma} \frac{|u(x)|^2}{(1+|x|)^2} \, dx
\leq C(b) \int_{\Gamma} \left( |u'(x)|^2 -
\lambda_b|u(x)|^2\right)\, dx, \quad u \in C_0^\infty(\Gamma_o),
\end{align*}
for a homogeneous tree $\Gamma$ as in the previous theorem.

\section{Proof of Theorems \ref{hardy_Mazya} and \ref{Thm:hardy_homo}}

\subsection{Orthogonal decomposition}

In this subsection we recall the results of Carlson \cite{C} and of Naimark and Solomyak \cite{NS2,NS1}. For each integer $k \geq 1$ we define the $k$-th branching functions $g_k:\R_+ \to \N$ by 
\begin{align*}
g_k(t) := \left\{
\begin{array}{l@{\quad}l}
 0,                         &          t < t_k     \, ,         \\
 1,                         & t_k \leq t < t_{k+1} \, ,         \\
 b_{k+1}b_{k+2} \cdots b_n, & t_n \leq t < t_{n+1} \, , k<n \, .
\end{array}
\right.
\end{align*}
The weighted Sobolev space $\H^1((t_k,\infty), g_k)$, $k\geq 0$, consists of all functions $f\in H^1_\loc(t_k,\infty)$ such that
\begin{align*}
\int_{t_k}^\infty \left( |f'(t)|^2 + |f(t)|^2 \right) g_k(t) \,dt  < \infty\, ,
\end{align*}
and we write $\H^1_0((t_k,\infty), g_k) := \{ f \in \H^1((t_k,\infty), g_k): f(t_k)=0\}$. Let $\A_k$ be the self-adjoint operator in $L_2((t_k,\infty),\, g_k)$
given by the closed quadratic form 
\begin{align*}
\qf_k[f] := \int_{t_k}^{\infty} |f'(t)|^2 g_k(t) \, dt
\end{align*}
with form domain $\H^1((0,\infty), g_0)$ if $k=0$ and form domain $\H^1_0((t_k,\infty), g_k)$ if $k \geq 1$. Notice that the operator $\A_0$ satisfies a \emph{Neumann} boundary condition at $t_0 = 0$, while the operators $\A_k$ with $k \geq 1$ satisfy \emph{Dirichlet} boundary conditions at $t_k$.

We call a function $V$ on $\Gamma$ \emph{symmetric} if it depends only on the distance from the root. In this case we shall sometimes abuse notation and write $V(|x|)$ instead of $V(x)$. The following statement is taken from \cite{NS1} and \cite{S}.

\begin{proposition} \label{NStheorem}
Let $\Gamma$ be a regular metric tree and let $V \in L_\infty(\Gamma)$ be symmetric. Then $-\Delta_\Neumann -
V$ is unitarily equivalent to the orthogonal sum of operators 
\begin{align} \label{decomp}
-\Delta_\Neumann - V \simeq (\A_0 - V) 
\oplus \sum_{k=1}^{\infty} 
\oplus \big(\A_k - V_k\big)^{[b_1...b_{k-1}(b_k-1)]}.
\end{align}
Here the symbol $[b_1...b_{k-1}(b_k-1)]$ means that the operator $\A_k
- V_k$ appears $b_1...b_{k-1}(b_k-1)$ times in the orthogonal sum, and
$V_k$ denotes the restriction of $V$ to the interval $(t_k,\infty)$.
\end{proposition}


\subsection{Proof of Theorem \ref{hardy_Mazya}}

Our proof follows closely the approach suggested in \cite{NS2} and is based on the following classical result by Muckenhoupt (see, e.g., \cite[Thm.~1.3.1/3]{M}).

\begin{proposition}\label{Mazya}
The inequality
\begin{align}\label{eq:Mazya}
\int_0^\infty \left|w(r)\int_r^\infty \varphi(s)\,ds\right|^2 \,dr
\leq S \int_0^\infty \left|v(r)\varphi(r)\right|^2 \,dr 
\end{align}
holds for all $\varphi$ if and only if
\begin{align}\label{eq:Mazya_crit}
T := \sup_{r>0} \int_0^r |w(s)|^2\, ds
\int_r^\infty \frac {ds}{|v(s)|^2} <\infty.
\end{align}
In this case, the sharp constant $S$ in \eqref{eq:Mazya} satisfies
\begin{align}\label{eq:Mazya_best}
T \leq S \leq 4 T.
\end{align}
\end{proposition}

\begin{proof}[Proof of Theorem \ref{hardy_Mazya}]
First we note that any function in $\H^1_0((t_k,\infty), g_k)$, $k\geq 1$, can be extended by zero to a function in $\H^1_0((0,\infty), g_0)$. In view of the definition of $g_k$ and the orthogonal decomposition from Proposition \ref{NStheorem} we see that inequality \eqref{eq:hardy_Mazya} is valid if and only if 
\begin{align}\label{eq:proof:hardy_Mazya1}
\int_{0}^\infty |f(t)|^2 \psi(t) g_0(t) \, dt \leq
C(\Gamma, \psi) \int_{0}^\infty |f'(t)|^2 g_0(t) \, dt,
\quad f \in \H^1((0,\infty), g_0)\, .
\end{align}
Since functions $f$ in $\H^1((0,\infty), g_0)$ can be represented as $f(t)=\int_t^\infty \varphi(s)\,ds$, the necessary and sufficient condition for inequality \eqref{eq:proof:hardy_Mazya1} follows from  Proposition \ref{Mazya} with $w = \sqrt{g_0 \psi}$ and $v = \sqrt{g_0}$. 
\end{proof}


\subsection{Proof of Theorem \ref{Thm:hardy_homo}}
Let $g_0$ be the first branching function of a homogeneous metric tree with edge length $1$ and branching number $b\geq 2$ and let $\lambda_b$ be the bottom of its
essential spectrum. Denote by $\omega$ the (unique up to a constant multiple) function on $\R_+$ satisfying
in distributional sense 
\begin{equation} \label{ground}
-\omega^{-1}(g_0 \omega')' = \lambda_b\, g_0 \, ,
\end{equation}
\begin{equation} \label{matching}
\omega'(0)=0, \quad \omega(j+)= \omega(j-),\quad
\omega'(j-)=b\, \omega'(j+), \quad j\in\N\, . 
\end{equation}
Using the explicit form of $\omega$ we showed in \cite{EFK} that $\sqrt{g_0} \omega$ grows linearly at infinity.

\begin{lemma} \label{efgrowth}
There exist constants $0<C_1<C_2<\infty$ (depending on $b$) such that
\begin{equation} \label{eq:efgrowth}
  C_1\, \frac{1+t} {\sqrt{g_0(t)}} \leq \omega(t)
  \leq C_2\, \frac{1+t}{\sqrt{g_0(t)}}\, , 
  \quad t\geq 0\, .
\end{equation}
\end{lemma}

Now we can prove a first version of Theorem \ref{Thm:hardy_homo}.

\begin{proposition} \label{hom}
The Hardy inequality
\begin{align}\label{eq:general_oneD_hardy}
\int_0^\infty \psi(t) |f(t)|^2 g_0(t) \, dt \leq C(g_0,\psi) \int_0^\infty
\left( |f'(t)|^2 - \lambda_b|f(t)|^2\right) g_0(t) \, dt\,,\\
f \in H^1(\R_+, g_0)\,, \notag
\end{align}
is valid if and only if
\begin{align*}
M'(g_0,\psi) := \sup_{r>0} \int_0^r \omega^2 \psi g_0 \, ds \int_r^\infty \frac
{ds}{\omega^2 g_0} < \infty.
\end{align*}
Moreover, the best constant in \eqref{eq:general_oneD_hardy} satisfies
\begin{align*}
M'(g_0,\psi) \leq C(g_0,\psi) \leq 4\,M'(g_0,\psi)\, .
\end{align*}
\end{proposition}

\begin{proof}
We know from Lemma \ref{efgrowth} that $\omega(t) >0$ for all
$t\in\R_+$. Therefore we can write any function $f \in H^1(\R_+, g_0)$
as a product $f= \omega\varphi$. Integrating by parts between the
points of discontinuity of $g_0$ and using \eqref{ground},
\eqref{matching} we arrive at the \emph{ground state representation}
\begin{equation}
\int_0^\infty
\left( |f'(t)|^2 - \lambda_b|f(t)|^2\right) g_0(t) \, dt =
\int_0^\infty\, |\omega(t)\varphi'(t)|^2g_0(t)\, dt\, .
\end{equation}
Similarly as in the proof of Theorem \ref{hardy_Mazya} we apply Proposition \ref{Mazya} with $v^2=\omega\, g_0$ and $w^2= \psi \omega^2\, g_0$ and obtain inequality \eqref{eq:general_oneD_hardy} with the claimed bounds on the constant.
\end{proof}

\begin{proof}[Proof of Theorem \ref{Thm:hardy_homo}]
By the same argument as in the proof of Theorem~\ref{hardy_Mazya}, inequality \eqref{eq:hardy_homo} is equivalent to the inequality
\begin{align*}
\int_0^\infty \psi(t) |f(t)|^2 g_0(t) \, dt \leq C(b,\psi) \int_0^\infty
\left( |f'(t)|^2 - \lambda_b|f(t)|^2\right) g_0(t) \, dt\,,
\quad f \in H^1(\R_+, g_0)\,.
\end{align*}
By Proposition \ref{hom} and Lemma \ref{efgrowth} the latter holds 
if and only if 
\begin{align*}
\sup_{r>0}\, (1+r)^{-1} \int_0^r \psi(t) (1+t)^2 \, dt<\infty\,,
\end{align*}
as claimed.
\end{proof}


\section{A loop graph with an Aharonov-Bohm magnetic field}

\subsection{Hardy's inequality on a loop graph}
Let $\Gamma=\Gamma_1\cup \Gamma_2$ be the graph embedded in $\R^2$ consisting of the circle
$$
\Gamma_1 :=\{ e^{i\theta}:\ 0\leq\theta\leq 2\pi\}
$$ 
around the origin of radius $1$ and the half-line
$$
\Gamma_2 = \{(r,0):\ r>1 \} \sim [1,\infty)\, .
$$
If $u$ is a function on $\Gamma$ we denote its restrictions to $\Gamma_1$ and $\Gamma_2$ by $u_1$ and $u_2$. For $a\in L_2(\Gamma_1)$ we consider the operator $H_a$ in $L_2(\Gamma)$ defined by the quadratic form
$$
h_a[u]:= \int_0^{2\pi}\, |-iu'_1(\theta) -a(\theta) u_1(\theta)|^2\, d\theta
+\int_1^\infty\, |u_2'(r)|^2\, dr\, .
$$
This operator describes a particle on $\Gamma$ subject to an Aharonov-Bohm magnetic field in the circle $\Gamma_1$ with flux 
\begin{equation}\label{eq:flux}
\alpha := \frac1{2\pi} \int_0^{2\pi} a(\theta)\,d\theta\,.
\end{equation}
It is easy to see that the Laplacian on $\Gamma$, i.e., $H_0$, does not satisfy a (non-trivial) Hardy inequality. We shall prove in this section that the situation is different for non-integer flux $\alpha$. To state our result we need

\begin{lemma}
	Let $0<\alpha\leq 1/2$. Then there is a unique solution $\lambda=\lambda_*(\alpha)$ of the equation
  \begin{equation} \label{implicit1}
		\cos(2\pi\alpha) =\cos(2\pi\sqrt{\lambda}) 
    - \frac{1-\sqrt{1-4\lambda}}{4\sqrt{\lambda}}\,
    \sin(2\pi\sqrt{\lambda})
  \end{equation}
  in the interval $(0,\alpha^2)$. 
\end{lemma}

It is easy to see that $\lambda_*(\alpha)$ is increasing with respect to $0<\alpha\leq 1/2$ and that $\lim_{\alpha\to 0} \lambda_*(\alpha)=0$. Numerically, one finds $\lambda_*(1/2)=0.1735$.

\begin{proof}
  Denoting the right hand side of \eqref{implicit1} by $h(\lambda)$, one checks that
  \begin{equation*}
    h(0)=1, \quad h(1/4)=-1,\quad h'(1/4)>0, \quad
    h''(\lambda)>0 \quad \forall\, \lambda\in(0,1/4)\, ,
  \end{equation*}
  which implies the assertion for $\alpha=1/2$. For $0<\alpha< 1/2$ one notes in addition that $h(0)>\cos(2\pi\alpha)>h(\alpha^2)$.
\end{proof}

In the following we extend $\lambda_*(\alpha)$ to an even, $1$-periodic function on $\R$. This function appears in the following Hardy inequality.

\begin{theorem} \label{AB1/2}
Let $a\in L_2(\Gamma_1)$ and define $\alpha$ by \eqref{eq:flux}. If $\alpha\not\in\Z$, then
\begin{align} \label{hardy-loop}
h_a[u] \geq \lambda_*(\alpha) \int_\Gamma \psi\, |u|^2 \,dx\,,
\quad u\in H^1(\Gamma)\,,
\end{align}
where $\psi(\theta):=1$ if $e^{i\theta}\in\Gamma_1$ and $\psi(r):=r^{-2}$ if $(r,0)\in\Gamma_2$. The constant $\lambda_*(\alpha)$ is sharp.
\end{theorem}

The proof of this inequality is relatively long and we break it into several steps.

\begin{figure}[htf] 
\includegraphics{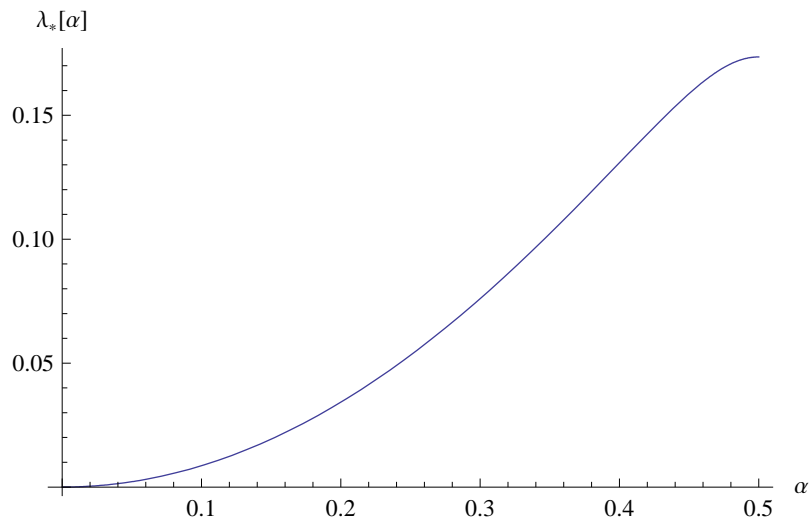}
\caption{The sharp constant $\lambda_*(\alpha)$ as function of $\alpha$}
\label{Fig.1}
\end{figure}


\subsection{The generalized groundstate}

In the following lemma we exhibit a function $\omega$ on $\Gamma$ which satisfies the Euler-Lagrange equation corresponding to the Hardy inequality \eqref{hardy-loop}.

\begin{lemma}\label{gs}
  Let $0<\alpha\leq 1/2$. There is a non-trivial function $\omega$ on $\Gamma$ such that
  \begin{align*}
  \left(-i\frac{d}{d\theta} -\alpha\right)^2\omega_1 = \lambda_*(\alpha) \, \omega_1
  & \quad \text{ on } \Gamma_1\,,\\
  -\frac{d^2}{dr^2}\omega_2 = \frac{\lambda_*(\alpha)}{r^2}\,\omega_2
  & \quad \text{ on } \Gamma_2
  \end{align*}
  and
  $$
  \omega_1(0)=\omega_1(2\pi)=\omega_2(1)\,,
  \quad \omega_1'(0)-\omega_1'(2\pi)+\omega_2'(1)=0\,.
  $$
  Moreover,
  \begin{enumerate}
  \item
  If $0<\alpha<1/2$, then $\omega(x)\neq 0$ for all $x\in\Gamma$, and if $\alpha=1/2$, then $\omega(x)= 0$ iff $x=(-1,0)\in\Gamma_1$.
  \item
  The function $\omega_2$ on $\Gamma_2$ is real-valued, and the function $j:=\re\overline{\omega_1}(-i\frac{d}{d\theta} -\alpha)\omega_1$ on $\Gamma_1$ is constant,
  $$
  j= - \sqrt{\lambda_*(\alpha)}\, \frac{\sin(2\pi\alpha)}{\sin(2\pi\sqrt{\lambda_*(\alpha)})}\,.
  $$
  \item
  If $0<\alpha<1/2$, then
  $$
  \int_0^{2\pi} \frac{d\theta}{|\omega_1|^2} 
  = \frac{2\pi\alpha}{\sqrt{\lambda_*(\alpha)}}\, \frac{\sin(2\pi\sqrt{\lambda_*(\alpha)})}{\sin(2\pi\alpha)}.
  $$
	\end{enumerate}
\end{lemma}

\begin{proof}
  Solving the equations separately on the circle $\Gamma_1$ and on the halfline $\Gamma_2$ and normalizing the functions to the value $1$ at the vertex we find that (with $\lambda_*(\alpha)=\mu^2$)
  \begin{align*}
    \omega_1(\theta) &
    = A e^{i\theta(\alpha+\mu)} + (1-A) e^{i\theta(\alpha-\mu)},
    \qquad A=\frac{e^{-2\pi i\alpha}-e^{-2\pi i\mu}}{2i\sin(2\pi\mu)}\,,\\
    \omega_2(r) &
    = r^\beta, \qquad  \beta =\frac{1-\sqrt{1-4\mu^2}}{2}\, .
  \end{align*}
  The matching condition for the derivatives is exactly the equation \eqref{implicit1} defining $\lambda_*(\alpha)$. Noting that
  $$
  |\omega_1(\theta)|^2 
  = \frac{2\,\sin(\pi(\alpha-\mu))\,\sin(\pi(\alpha+\mu))}{\sin^2(2\pi\mu)} \left(B-\cos(2\mu(\theta-\pi))\right)
  $$
  where
  $$
  B:=\frac{\sin^2(\pi(\alpha-\mu)) +
  \sin^2(\pi(\alpha+\mu))}{2\,\sin(\pi(\alpha-\mu))\,\sin(\pi(\alpha+\mu))}\,,
  $$
  we see easily that $\omega_1(\theta)$ vanishes only when $\alpha=1/2$ and $\theta=\pi$. Moreover, for $0<\alpha<1/2$ we find
  \begin{align*}
  \int_0^{2\pi} \frac{d\theta}{|\omega_1|^2} 
  & = \frac{\sin^2(2\pi\mu)}{\sin(\pi(\alpha-\mu))\,\sin(\pi(\alpha+\mu))}
	\frac{1}{\mu \sqrt{B^2-1}} \arctan \left(\tan(\mu\pi) \sqrt{\frac{B+1}{B-1}} \right) \\
	& = 2\pi\alpha\, \frac{\sin(2\pi\mu)}{\mu \sin(2\pi\alpha)}\,.
  \end{align*}
  In the last step we used trigonometric identities to simplify the expressions for $B\pm1$. It is elementary to check that $\re\overline{\omega_1}(-i\frac{d}{d\theta} -\alpha)\omega_1$ is constant and has the value given in the lemma.
\end{proof}

Next we derive a ground state representation formula. Note that an additional term appears in the magnetic case.

\begin{lemma}\label{gsr}
	Let $0<\alpha\leq 1/2$ and let $\omega$ and $j$ be as in Lemma \ref{gs}. Let $u\in H^1(\Gamma)$ and, in case $\alpha=1/2$, assume that $u_1(\pi)=0$. Then with $v:=\omega^{-1} u$,
	\begin{equation}\label{eq:gsr}
	h_\alpha[u] - \lambda_*(\alpha) \int_\Gamma \psi\, |u|^2 \,dx 
	= \int_\Gamma |\omega|^2 |v'|^2\,dx 
	+ 2 j\im \int_0^{2\pi} \overline{v_1}v_1'\, d\theta\,. 
	\end{equation}
\end{lemma}

This follows using integration by parts and the equation satisfied by $\omega$. We omit the details. In order to deal with the second term on the right hand side of \eqref{eq:gsr} we shall need

\begin{lemma}\label{change}
Let $w$ be a non-negative function on $(0,2\pi)$.
Then for any periodic $H^1$-function $v$ on $(0,2\pi)$,
\begin{equation}\label{eq:change}
\int_0^{2\pi} |v'|^2 w \,d\theta
\geq 2\pi \left( \int_0^{2\pi} \frac{d\theta}{w} \right)^{-1}
\left|\im \int_0^{2\pi} \overline v v'\,d\theta\right| \, .
\end{equation}
\end{lemma}

\begin{proof}
Let $\Phi(\theta):= \int_0^\theta
w(\theta)^{-1}\,d\theta$ and $l:=\Phi(2\pi)$. Then $\Phi$
is a strictly increasing function which maps $(0,2\pi)$ onto
  $(0,l)$. Given a periodic $H^1$-function $v$ on $(0,2\pi)$,
        define $f$ on $(0,l)$ by $f(\Phi(\theta)):=v(\theta)$. Then
        $f$ is a periodic $H^1$-function $v$ on $(0,l)$ and one has 
	$$
		v'(\theta)=f'(\Phi(\theta)) \Phi'(\theta) =
                f'(\Phi(\theta)) w(\theta)^{-1}\,. 
	$$
	Hence 
	$$
		\int_0^{2\pi} |v'|^2\,w\,d\theta = \int_0^l |f'|^2\,d\phi\,,
		\quad
		\int_0^{2\pi} \overline v v' \,d\theta = \int_0^l \overline ff'\,d\phi\,.
	$$
	and for any $\beta\in\R$
	\begin{align*}
		& \int_0^{2\pi} |v'|^2\,w\,d\theta 
		- 4\pi\beta l^{-1} \im \int_0^{2\pi} \overline vv'\,d\theta \\
		& \quad = \int_0^l \left|f'- i\frac{2\pi\beta}{l} f\right|^2\,d\phi 
		- 4\pi^2\left(\frac\beta l \right)^2 \int_0^l |f|^2\,d\phi \\
		& \quad = 4\pi^2l^{-2} \sum_{n\in\Z} \left( |n-\beta|^2 -
                  |\beta|^2 \right) |\hat f_n|^2\,, 
	\end{align*}
	where $\hat f_n:=l^{-1/2} \int_0^l f(\phi) e^{-i2\pi n\phi/l}\,d\phi$ are the Fourier coefficients of $f$. Clearly,      this is non-negative for $|\beta|\leq 1/2$. Choosing $\beta=\pm 1/2$ we obtain the assertion.
\end{proof}


\subsection{Proof of Theorem \ref{AB1/2}}

\emph{First step.} We claim that it is enough to consider the case $a\equiv\alpha\in(0,1/2]$. Indeed, since multiplication by $\exp(-i\alpha\theta+i\int_0^\theta a(\theta')\,d\theta')$ on $\Gamma_1$ is unitary, we see that $H_a$ and $H_\alpha$ are unitarily equivalent. Similarly, since multiplication on $\Gamma_1$ by $e^{in\theta}$ is unitary, $H_\alpha$ and $H_{\alpha'}$ are unitarily equivalent if $\alpha-\alpha'\in\Z$. Moreover, $H_\alpha$ and $H_{-\alpha}$ are anti-unitarily equivalent by complex conjugation. This proves the claim.

\emph{Second step.} We prove the assertion for $0<\alpha<1/2$. Let $\omega$ and $j$ be as in Lemma~\ref{gs}. Then Lemmas \ref{gsr} and \ref{change} imply that for any $u=v\omega\in H^1(\Gamma)$,
$$
h_\alpha[u] - \lambda_*(\alpha) \int_\Gamma \psi\, |u|^2 \,dx 
	\geq \left(1- \frac{|j|}\pi \int_0^{2\pi} \frac{d\theta}{|\omega|^2} \right)
	 \int_\Gamma |\omega|^2 |v'|^2\,dx \,. 
$$	
This is non-negative in view of the explicit expressions in Lemma \ref{gs}, proving \eqref{hardy-loop}. Letting $v$ approach the constant function, we find that the inequality is sharp. Note that the previous argument does not work for $\alpha=1/2$ since $\omega_1(\pi)=0$ in this case.

	\emph{Third step.} In the remainder of the proof we assume that $\alpha=1/2$. We show that \eqref{hardy-loop} is equivalent to two independent inequalities for functions satisfying appropriate Dirichlet boundary conditions. For this purpose, we decompose $u=u^s+u^a$ into its (twisted) symmetric and antisymmetric parts,
	\begin{align*}
	u_1^s := \frac12 (u_1(\theta)+e^{i\theta} u_1(2\pi-\theta))\,,
	\quad & u_2^s(r):=u_2(r)\,, \\
	u_1^a := \frac12 (u_1(\theta)-e^{i\theta} u_1(2\pi-\theta))\,,
	\quad & u_2^a(r):=0\,.
	\end{align*}
	An easy calculation shows that
	$$
  h_{1/2}[u] = h_{1/2}[u^s] + h_{1/2}[u^a]\,,
  \quad
  \int_\Gamma \psi\, |u|^2 \,dx 
  = \int_\Gamma \psi\, |u^a|^2 \,dx + \int_\Gamma \psi\, |u^s|^2 \,dx\, .
  $$
	Hence the desired inequality decouples into two independent inequalities,
	\begin{align}
		\label{eq:absymm}
		h_{1/2}[u^s] & \geq \lambda_*(1/2) \int_\Gamma \psi\, |u^s|^2 \,dx\,, 
		\quad u^s_1(\pi)=0\,,\\
		\label{eq:abantisymm}
		\int_0^{2\pi}\, \left|\left(-i\frac{d}{d\theta}-\frac 12\right) u_1^a\right|^2\, d\theta
		& \geq \lambda_*(1/2) \int_0^{2\pi}\, |u_1^a|^2\, d\theta\,,
		\quad u_1^a(0)=u_1^a(2\pi)=0.
	\end{align}
	To prove \eqref{eq:abantisymm} we write $u_1^a(\theta)=e^{i\theta/2} f(\theta)$. Using that the first Dirichlet eigenvalue on $(0,2\pi)$ is $1/4$ we find
	$$
		\int_0^{2\pi}\, \left|\left(-i\frac{d}{d\theta}-\frac 12\right) u_1^a\right|^2\, d\theta
		= \int_0^{2\pi}\, |f'|^2\, d\theta \geq \frac 14 \int_0^{2\pi}\, |f|^2\, d\theta
		= \frac 14 \int_0^{2\pi}\, |u_1^a|^2\, d\theta\,.
	$$
	Since $\lambda_*(1/2)<1/4$, we obtain \eqref{eq:abantisymm}.
	
	To prove \eqref{eq:absymm} we argue similarly as in the second step. Let again $\omega$ and $j$ be as in Lemma \ref{gs}. Since $u_1^s$ vanishes at $\theta=\pi$ and $\omega$ is non-zero away from this point, we can write $u^s = v\omega$ and apply Lemma \ref{gsr}. Note that $j=0$ for $\alpha=1/2$, which proves \eqref{eq:absymm}. That the constant $\lambda_*(1/2)$ is sharp follows again by letting $v$ approach the constant function. This concludes the proof of the theorem.

\bibliographystyle{amsalpha}

\begin{thebibliography}{GGMT}

\bibitem[C]{C}
R.~Carlson, Nonclassical Strum-Liouville problems and Schr\"odinger operators
on radial trees. {\em Electron J. Differential Equation} {\bf 71}
(2000), 24pp.

\bibitem[EFK]{EFK} T.~Ekholm, R.~L.~Frank, H.~Kova\v r\'{\i}k,
Eigenvalue estimates for Schr\"odinger operators on metric
trees. Preprint: arXiv: math.SP/0710.5500v1.  

\bibitem[EHP]{EHP} W.~D.~Evans, D.~J.~Harris and L.~Pick, Weighted Hardy
  and Poincar\'e inequalities on trees. {\em J.~London Math.~Soc.} (2)
  {\bf 52} (1995), no.1, 121--136.

\bibitem[K]{K} H.~Kova\v r\'{\i}k, Weakly coupled Schr\"odinger operators
  on regular metric trees. Preprint: arXiv: math-ph/0608013.

\bibitem[LW]{LW}
	A. Laptev and T. Weidl, Hardy inequalities for magnetic Dirichlet forms, \textit{Operator Theory: Advances and Applications} \textbf{108}, Birkh\"auser, Basel (1999), 299--305,

\bibitem[M]{M}
V.~G.~Mazya, {\em Sobolev spaces}, Springer, Berlin, Heidelberg (1985).

\bibitem[NS1]{NS2}
K.~Naimark and M.~Solomyak, Eigenvalue estimates for the weighted
Laplacian on metric trees, {\em Proc. London  Math. Soc.} {\bf 80}
(2000), no. 3, 690--724.

\bibitem[NS2]{NS1}
K.~Naimark and M.~Solomyak, Geometry of the Sobolev spaces on the
regular trees and Hardy's inequalities, {\em Russ. J. Math. Phys.}
{\bf 8} (2001), no. 3, 322--335.


\bibitem[SS]{SS}
A.~Sobolev and M.~Solomyak, Schr\"odinger operators on homogeneous
metric trees: spectrum in gaps,
{\em Rev. Math. Phys.} {\bf 14} (2002), no. 5, 421--467.

\bibitem[S]{S}
M.~Solomyak, On the spectrum of the Laplacian on metric trees. Special
section on quantum graphs. {\em Waves Random Media} {\bf 14} (2004),
no. 1, S155--S171.

\end{thebibliography}

\end{document}